\documentclass{arxiv_elsart}

\usepackage{graphicx}

\begin{document}

\begin{frontmatter}

\title{Note on the Irreducible Triangulations of the Klein Bottle}

\author{Thom Sulanke}
\address{Department of Physics, Indiana University, Bloomington, 
	Indiana 47405}
\ead{tsulanke@indiana.edu}

\begin{abstract}
We give the complete list of the $29$ irreducible triangulations 
of the Klein bottle.
We show how the construction of Lawrencenko and Negami, which
listed only $25$ such irreducible triangulations, can
be modified at two points to produce the $4$ additional irreducible 
triangulations of the Klein bottle.

\end{abstract}

\end{frontmatter}

\section{Introduction}
\label{intro}

A {\em triangulation\/} of a closed surface 
is a simple graph embedded in the surface
so that each face is a triangle and so that any two faces share at most
one edge.   Two triangulations $G$ and $G'$ of
a surface are {\em equivalent\/} 
if there is a homeomorphism $h$
with $h(G)=G'$.  

Let $e = ac$ be an edge in a triangulation $G$ and
$abc$ and $acd$ be the two faces which have $e$ as a common edge.
The {\em contraction\/} of $e$ is obtained by deleting $ac$,
identifying vertices $a$ and $c$, 
removing one of the multiple edges $ab$ or $cb$,
and removing one of the multiple edges $ad$ or $cd$.
An edge $e$ of a triangulation $G$ is {\em contractible\/} if the 
contraction of $e$ yields another triangulation of the surface in
which $G$ is embedded.  
If an edge $e$ is contained in a three cycle other than the two 
which bound the faces which share $e$, then the contraction of $e$ 
produces multiple edges.
So, for a triangulation $G$, not $K_4$ embedded in the sphere, 
an edge $e$ of $G$ 
is not contractible
if and only if $e$ is contained in at least three cycles of length $3$. 
A triangulation is said to be {\em irreducible\/} 
if it has no contractible edge. 

This note assumes the reader's familiarity with Lawrencenko and Negami's
paper~\cite{LN} where their Theorem 1 claims the existence of exactly 25
irreducible triangulations of the Klein bottle.
Theorem~\ref{thm1} corrects this claim; it requires
Lemma~\ref{lem6}, a modification of Lemma 6 of~\cite{LN},
which claimed the existence of exactly 21 irreducible 
triangulations of handle type.
All other results of~\cite{LN} remain valid.

\begin{thm}
\label{thm1}
There are exactly 29 nonequivalent irreducible triangulations of the Klein
bottle.
\end{thm}

\begin{lem}
\label{lem6}
There are exactly 25 nonequivalent irreducible triangulations of handle
type.
\end{lem}

In Section~\ref{list}
we list all 29 irreducible triangulations of the Klein bottle
and show that each of the new ones is not equivalent 
to any of the others.
A complete proof that these 29 irreducible triangulations
are all of the irreducible triangulations of the Klein bottle requires 
repeating the construction in \cite{LN} along with two modifications.  
We give only the modifications required.
In Section~\ref{mod-Kh25} we examine a subcase which was overlooked 
in~\cite{LN} and which leads to the additional triangulation Kh25.
In Section~\ref{mod-Kh222324} we examine in more detail than was 
given in~\cite{LN} a subcase which leads to the additional triangulations 
Kh22, Kh23, and Kh24.

\section{List of Irreducible Triangulations}
\label{list}

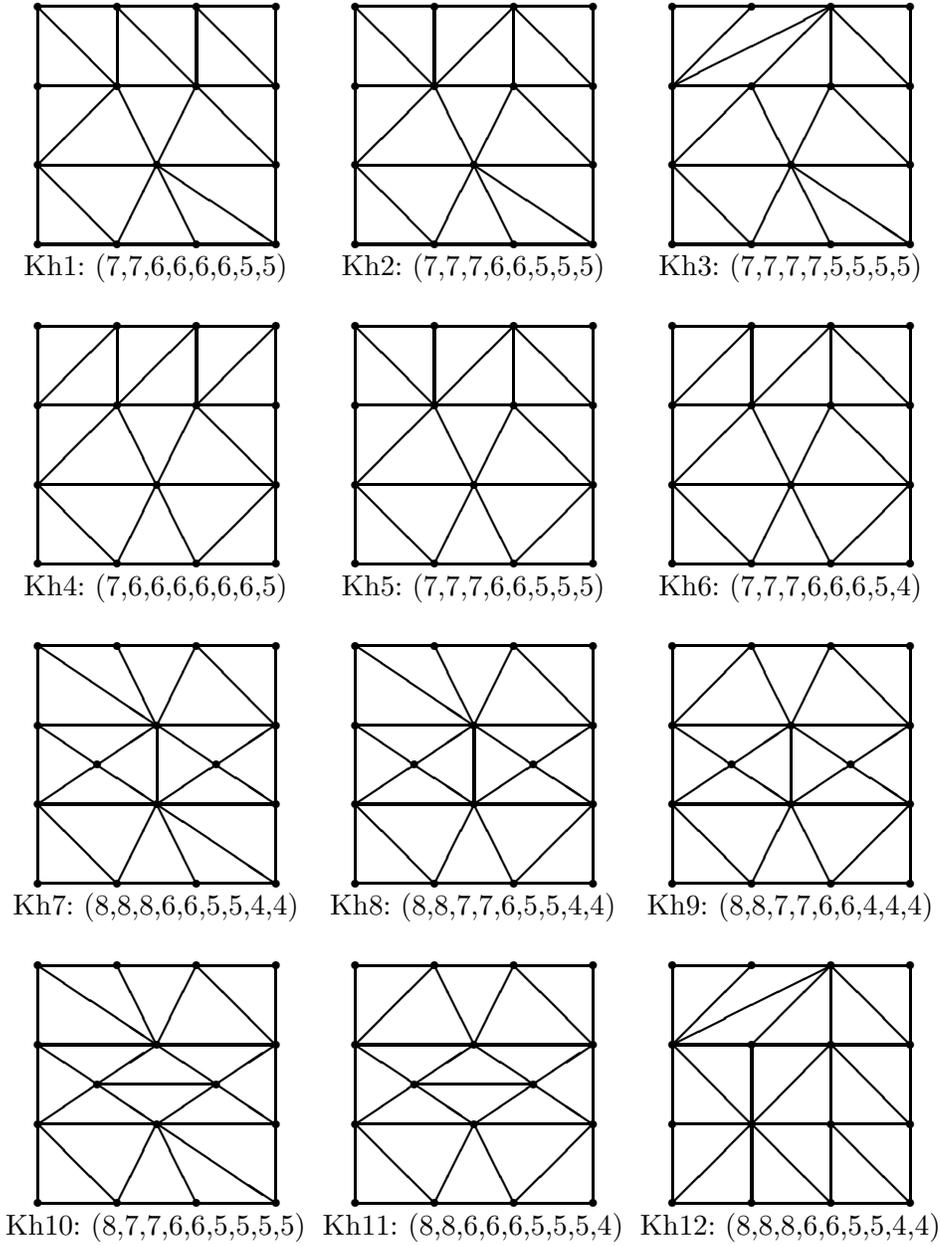
\begin{figure}
\begin{center}
\setlength{\unitlength}{3pt}
\small

\begin{picture}(110,40)(0,0)

\put(0,0){
\begin{picture}(30,40)(0,-5)
\thicklines
\multiput(0,0)(0,10){4}{\line(1,0){30}}
\put(0,0){\line(0,1){30}}
\put(30,0){\line(0,1){30}}
\multiput(0,0)(10,0){4}{\circle*{1}}
\multiput(0,10)(15,0){3}{\circle*{1}}
\multiput(0,20)(10,0){4}{\circle*{1}}
\multiput(0,30)(10,0){4}{\circle*{1}}
\put(0,10){\line(1,-1){10}}
\put(10,0){\line(1,2){5}}
\put(15,10){\line(1,-2){5}}
\put(15,10){\line(3,-2){15}}
\put(0,10){\line(1,1){10}}
\put(10,20){\line(1,-2){5}}
\put(15,10){\line(1,2){5}}
\put(20,20){\line(1,-1){10}}
\put(0,30){\line(1,-1){10}}
\put(10,20){\line(0,1){10}}
\put(10,30){\line(1,-1){10}}
\put(20,20){\line(0,1){10}}
\put(20,30){\line(1,-1){10}}
\put(15,-3){\makebox(0,0)[]{Kh1: (7,7,6,6,6,6,5,5)}}
\end{picture}
}

\put(40,0){
\begin{picture}(30,40)(0,-5)
\thicklines
\multiput(0,0)(0,10){4}{\line(1,0){30}}
\put(0,0){\line(0,1){30}}
\put(30,0){\line(0,1){30}}
\multiput(0,0)(10,0){4}{\circle*{1}}
\multiput(0,10)(15,0){3}{\circle*{1}}
\multiput(0,20)(10,0){4}{\circle*{1}}
\multiput(0,30)(10,0){4}{\circle*{1}}
\put(0,10){\line(1,-1){10}}
\put(10,0){\line(1,2){5}}
\put(15,10){\line(1,-2){5}}
\put(15,10){\line(3,-2){15}}
\put(0,10){\line(1,1){10}}
\put(10,20){\line(1,-2){5}}
\put(15,10){\line(1,2){5}}
\put(20,20){\line(1,-1){10}}
\put(0,30){\line(1,-1){10}}
\put(10,20){\line(0,1){10}}
\put(10,20){\line(1,1){10}}
\put(20,20){\line(0,1){10}}
\put(20,30){\line(1,-1){10}}
\put(15,-3){\makebox(0,0)[]{Kh2: (7,7,7,6,6,5,5,5)}}
\end{picture}
}

\put(80,0){
\begin{picture}(30,40)(0,-5)
\thicklines
\multiput(0,0)(0,10){4}{\line(1,0){30}}
\put(0,0){\line(0,1){30}}
\put(30,0){\line(0,1){30}}
\multiput(0,0)(10,0){4}{\circle*{1}}
\multiput(0,10)(15,0){3}{\circle*{1}}
\multiput(0,20)(10,0){4}{\circle*{1}}
\multiput(0,30)(10,0){4}{\circle*{1}}
\put(0,10){\line(1,-1){10}}
\put(10,0){\line(1,2){5}}
\put(15,10){\line(1,-2){5}}
\put(15,10){\line(3,-2){15}}
\put(0,10){\line(1,1){10}}
\put(10,20){\line(1,-2){5}}
\put(15,10){\line(1,2){5}}
\put(20,20){\line(1,-1){10}}
\put(0,20){\line(1,1){10}}
\put(0,20){\line(2,1){20}}
\put(10,20){\line(1,1){10}}
\put(20,20){\line(0,1){10}}
\put(20,30){\line(1,-1){10}}
\put(15,-3){\makebox(0,0)[]{Kh3: (7,7,7,7,5,5,5,5)}}
\end{picture}
}

\end{picture}

\begin{picture}(110,40)(0,0)

\put(0,0){
\begin{picture}(30,40)(0,-5)
\thicklines
\multiput(0,0)(0,10){4}{\line(1,0){30}}
\put(0,0){\line(0,1){30}}
\put(30,0){\line(0,1){30}}
\multiput(0,0)(10,0){4}{\circle*{1}}
\multiput(0,10)(15,0){3}{\circle*{1}}
\multiput(0,20)(10,0){4}{\circle*{1}}
\multiput(0,30)(10,0){4}{\circle*{1}}
\put(0,10){\line(1,-1){10}}
\put(10,0){\line(1,2){5}}
\put(15,10){\line(1,-2){5}}
\put(20,0){\line(1,1){10}}
\put(0,10){\line(1,1){10}}
\put(10,20){\line(1,-2){5}}
\put(15,10){\line(1,2){5}}
\put(20,20){\line(1,-1){10}}
\put(0,20){\line(1,1){10}}
\put(10,20){\line(0,1){10}}
\put(10,20){\line(1,1){10}}
\put(20,20){\line(0,1){10}}
\put(20,20){\line(1,1){10}}
\put(15,-3){\makebox(0,0)[]{Kh4: (7,6,6,6,6,6,6,5)}}
\end{picture}
}

\put(40,0){
\begin{picture}(30,40)(0,-5)
\thicklines
\multiput(0,0)(0,10){4}{\line(1,0){30}}
\put(0,0){\line(0,1){30}}
\put(30,0){\line(0,1){30}}
\multiput(0,0)(10,0){4}{\circle*{1}}
\multiput(0,10)(15,0){3}{\circle*{1}}
\multiput(0,20)(10,0){4}{\circle*{1}}
\multiput(0,30)(10,0){4}{\circle*{1}}
\put(0,10){\line(1,-1){10}}
\put(10,0){\line(1,2){5}}
\put(15,10){\line(1,-2){5}}
\put(20,0){\line(1,1){10}}
\put(0,10){\line(1,1){10}}
\put(10,20){\line(1,-2){5}}
\put(15,10){\line(1,2){5}}
\put(20,20){\line(1,-1){10}}
\put(0,30){\line(1,-1){10}}
\put(10,20){\line(0,1){10}}
\put(10,20){\line(1,1){10}}
\put(20,20){\line(0,1){10}}
\put(20,30){\line(1,-1){10}}
\put(15,-3){\makebox(0,0)[]{Kh5: (7,7,7,6,6,5,5,5)}}
\end{picture}
}

\put(80,0){
\begin{picture}(30,40)(0,-5)
\thicklines
\multiput(0,0)(0,10){4}{\line(1,0){30}}
\put(0,0){\line(0,1){30}}
\put(30,0){\line(0,1){30}}
\multiput(0,0)(10,0){4}{\circle*{1}}
\multiput(0,10)(15,0){3}{\circle*{1}}
\multiput(0,20)(10,0){4}{\circle*{1}}
\multiput(0,30)(10,0){4}{\circle*{1}}
\put(0,10){\line(1,-1){10}}
\put(10,0){\line(1,2){5}}
\put(15,10){\line(1,-2){5}}
\put(20,0){\line(1,1){10}}
\put(0,10){\line(1,1){10}}
\put(10,20){\line(1,-2){5}}
\put(15,10){\line(1,2){5}}
\put(20,20){\line(1,-1){10}}
\put(0,20){\line(1,1){10}}
\put(10,20){\line(0,1){10}}
\put(10,20){\line(1,1){10}}
\put(20,20){\line(0,1){10}}
\put(20,30){\line(1,-1){10}}
\put(15,-3){\makebox(0,0)[]{Kh6: (7,7,7,6,6,6,5,4)}}
\end{picture}
}

\end{picture}

\begin{picture}(110,40)(0,0)

\put(0,0){
\begin{picture}(30,40)(0,-5)
\thicklines
\multiput(0,0)(0,10){4}{\line(1,0){30}}
\put(0,0){\line(0,1){30}}
\put(30,0){\line(0,1){30}}
\multiput(0,0)(10,0){4}{\circle*{1}}
\multiput(0,10)(15,0){3}{\circle*{1}}
\multiput(7.5,15)(15,0){2}{\circle*{1}}
\multiput(0,20)(15,0){3}{\circle*{1}}
\multiput(0,30)(10,0){4}{\circle*{1}}
\put(0,10){\line(1,-1){10}}
\put(10,0){\line(1,2){5}}
\put(15,10){\line(1,-2){5}}
\put(15,10){\line(3,-2){15}}
\put(0,10){\line(3,2){15}}
\put(0,20){\line(3,-2){15}}
\put(15,10){\line(0,1){10}}
\put(15,10){\line(3,2){15}}
\put(15,20){\line(3,-2){15}}
\put(0,30){\line(3,-2){15}}
\put(10,30){\line(1,-2){5}}
\put(15,20){\line(1,2){5}}
\put(20,30){\line(1,-1){10}}
\put(15,-3){\makebox(0,0)[]{Kh7: (8,8,8,6,6,5,5,4,4)}}
\end{picture}
}

\put(40,0){
\begin{picture}(30,40)(0,-5)
\thicklines
\multiput(0,0)(0,10){4}{\line(1,0){30}}
\put(0,0){\line(0,1){30}}
\put(30,0){\line(0,1){30}}
\multiput(0,0)(10,0){4}{\circle*{1}}
\multiput(0,10)(15,0){3}{\circle*{1}}
\multiput(7.5,15)(15,0){2}{\circle*{1}}
\multiput(0,20)(15,0){3}{\circle*{1}}
\multiput(0,30)(10,0){4}{\circle*{1}}
\put(0,10){\line(1,-1){10}}
\put(10,0){\line(1,2){5}}
\put(15,10){\line(1,-2){5}}
\put(20,0){\line(1,1){10}}
\put(0,10){\line(3,2){15}}
\put(0,20){\line(3,-2){15}}
\put(15,10){\line(0,1){10}}
\put(15,10){\line(3,2){15}}
\put(15,20){\line(3,-2){15}}
\put(0,30){\line(3,-2){15}}
\put(10,30){\line(1,-2){5}}
\put(15,20){\line(1,2){5}}
\put(20,30){\line(1,-1){10}}
\put(15,-3){\makebox(0,0)[]{Kh8: (8,8,7,7,6,5,5,4,4)}}
\end{picture}
}

\put(80,0){
\begin{picture}(30,40)(0,-5)
\thicklines
\multiput(0,0)(0,10){4}{\line(1,0){30}}
\put(0,0){\line(0,1){30}}
\put(30,0){\line(0,1){30}}
\multiput(0,0)(10,0){4}{\circle*{1}}
\multiput(0,10)(15,0){3}{\circle*{1}}
\multiput(7.5,15)(15,0){2}{\circle*{1}}
\multiput(0,20)(15,0){3}{\circle*{1}}
\multiput(0,30)(10,0){4}{\circle*{1}}
\put(0,10){\line(1,-1){10}}
\put(10,0){\line(1,2){5}}
\put(15,10){\line(1,-2){5}}
\put(20,0){\line(1,1){10}}
\put(0,10){\line(3,2){15}}
\put(0,20){\line(3,-2){15}}
\put(15,10){\line(0,1){10}}
\put(15,10){\line(3,2){15}}
\put(15,20){\line(3,-2){15}}
\put(0,20){\line(1,1){10}}
\put(10,30){\line(1,-2){5}}
\put(15,20){\line(1,2){5}}
\put(20,30){\line(1,-1){10}}
\put(15,-3){\makebox(0,0)[]{Kh9: (8,8,7,7,6,6,4,4,4)}}
\end{picture}
}

\end{picture}

\begin{picture}(110,40)(0,0)

\put(0,0){
\begin{picture}(30,40)(0,-5)
\thicklines
\multiput(0,0)(0,10){4}{\line(1,0){30}}
\put(0,0){\line(0,1){30}}
\put(30,0){\line(0,1){30}}
\multiput(0,0)(10,0){4}{\circle*{1}}
\multiput(0,10)(15,0){3}{\circle*{1}}
\multiput(7.5,15)(15,0){2}{\circle*{1}}
\multiput(0,20)(15,0){3}{\circle*{1}}
\multiput(0,30)(10,0){4}{\circle*{1}}
\put(0,10){\line(1,-1){10}}
\put(10,0){\line(1,2){5}}
\put(15,10){\line(1,-2){5}}
\put(15,10){\line(3,-2){15}}
\put(0,10){\line(3,2){15}}
\put(0,20){\line(3,-2){15}}
\put(7.5,15){\line(1,0){15}}
\put(15,10){\line(3,2){15}}
\put(15,20){\line(3,-2){15}}
\put(0,30){\line(3,-2){15}}
\put(10,30){\line(1,-2){5}}
\put(15,20){\line(1,2){5}}
\put(20,30){\line(1,-1){10}}
\put(15,-3){\makebox(0,0)[]{Kh10: (8,7,7,6,6,5,5,5,5)}}
\end{picture}
}

\put(40,0){
\begin{picture}(30,40)(0,-5)
\thicklines
\multiput(0,0)(0,10){4}{\line(1,0){30}}
\put(0,0){\line(0,1){30}}
\put(30,0){\line(0,1){30}}
\multiput(0,0)(10,0){4}{\circle*{1}}
\multiput(0,10)(15,0){3}{\circle*{1}}
\multiput(7.5,15)(15,0){2}{\circle*{1}}
\multiput(0,20)(15,0){3}{\circle*{1}}
\multiput(0,30)(10,0){4}{\circle*{1}}
\put(0,10){\line(1,-1){10}}
\put(10,0){\line(1,2){5}}
\put(15,10){\line(1,-2){5}}
\put(20,0){\line(1,1){10}}
\put(0,10){\line(3,2){15}}
\put(0,20){\line(3,-2){15}}
\put(7.5,15){\line(1,0){15}}
\put(15,10){\line(3,2){15}}
\put(15,20){\line(3,-2){15}}
\put(0,20){\line(1,1){10}}
\put(10,30){\line(1,-2){5}}
\put(15,20){\line(1,2){5}}
\put(20,30){\line(1,-1){10}}
\put(15,-3){\makebox(0,0)[]{Kh11: (8,8,6,6,6,5,5,5,4)}}
\end{picture}
}

\put(80,0){
\begin{picture}(30,40)(0,-5)
\thicklines
\multiput(0,0)(0,10){4}{\line(1,0){30}}
\put(0,0){\line(0,1){30}}
\put(30,0){\line(0,1){30}}
\multiput(0,0)(10,0){4}{\circle*{1}}
\multiput(0,10)(10,0){4}{\circle*{1}}
\multiput(0,20)(10,0){4}{\circle*{1}}
\multiput(0,30)(10,0){4}{\circle*{1}}
\put(0,0){\line(1,1){10}}
\put(10,0){\line(0,1){10}}
\put(10,10){\line(1,-1){10}}
\put(20,0){\line(0,1){10}}
\put(20,10){\line(1,-1){10}}
\put(0,20){\line(1,-1){10}}
\put(10,10){\line(0,1){10}}
\put(10,10){\line(1,1){10}}
\put(20,10){\line(0,1){10}}
\put(20,20){\line(1,-1){10}}
\put(0,20){\line(1,1){10}}
\put(0,20){\line(2,1){20}}
\put(10,20){\line(1,1){10}}
\put(20,20){\line(0,1){10}}
\put(20,30){\line(1,-1){10}}
\put(15,-3){\makebox(0,0)[]{Kh12: (8,8,8,6,6,5,5,4,4)}}
\end{picture}
}

\end{picture}

\caption{Irreducible triangulations of the Klein bottle, Kh1 - Kh12}
\label{fig1}
\end{center}
\end{figure}

\begin{figure}
\begin{center}
\setlength{\unitlength}{3pt}
\small

\begin{picture}(110,40)(0,0)

\put(0,0){
\begin{picture}(30,40)(0,-5)
\thicklines
\multiput(0,0)(0,10){4}{\line(1,0){30}}
\put(0,0){\line(0,1){30}}
\put(30,0){\line(0,1){30}}
\multiput(0,0)(10,0){4}{\circle*{1}}
\multiput(0,10)(10,0){4}{\circle*{1}}
\multiput(0,20)(10,0){4}{\circle*{1}}
\multiput(0,30)(10,0){4}{\circle*{1}}
\put(0,0){\line(1,1){10}}
\put(10,0){\line(0,1){10}}
\put(10,10){\line(1,-1){10}}
\put(20,0){\line(0,1){10}}
\put(20,10){\line(1,-1){10}}
\put(0,20){\line(1,-1){10}}
\put(10,10){\line(0,1){10}}
\put(10,10){\line(1,1){10}}
\put(20,10){\line(0,1){10}}
\put(20,20){\line(1,-1){10}}
\put(0,20){\line(1,1){10}}
\put(0,20){\line(2,1){20}}
\put(10,20){\line(1,1){10}}
\put(20,20){\line(0,1){10}}
\put(20,20){\line(1,1){10}}
\put(15,-3){\makebox(0,0)[]{Kh13: (8,8,7,7,7,5,4,4,4)}}
\end{picture}
}

\put(40,0){
\begin{picture}(30,40)(0,-5)
\thicklines
\multiput(0,0)(0,10){4}{\line(1,0){30}}
\multiput(0,0)(10,0){4}{\line(0,1){30}}
\multiput(0,0)(10,0){4}{\circle*{1}}
\multiput(0,10)(10,0){4}{\circle*{1}}
\multiput(0,20)(10,0){4}{\circle*{1}}
\multiput(0,30)(10,0){4}{\circle*{1}}
\put(0,0){\line(1,1){10}}
\put(10,0){\line(1,1){10}}
\put(20,0){\line(1,1){10}}
\put(0,10){\line(1,1){10}}
\put(10,10){\line(1,1){10}}
\put(20,10){\line(1,1){10}}
\put(0,20){\line(1,1){10}}
\put(10,20){\line(1,1){10}}
\put(20,20){\line(1,1){10}}
\put(15,-3){\makebox(0,0)[]{Kh14: (6,6,6,6,6,6,6,6,6)}}
\end{picture}
}

\put(80,0){
\begin{picture}(30,40)(0,-5)
\thicklines
\multiput(0,0)(0,10){4}{\line(1,0){30}}
\multiput(0,0)(10,0){4}{\line(0,1){30}}
\multiput(0,0)(10,0){4}{\circle*{1}}
\multiput(0,10)(10,0){4}{\circle*{1}}
\multiput(0,20)(10,0){4}{\circle*{1}}
\multiput(0,30)(10,0){4}{\circle*{1}}
\put(0,10){\line(1,-1){10}}
\put(10,0){\line(1,1){10}}
\put(20,0){\line(1,1){10}}
\put(0,10){\line(1,1){10}}
\put(10,10){\line(1,1){10}}
\put(20,10){\line(1,1){10}}
\put(0,20){\line(1,1){10}}
\put(10,20){\line(1,1){10}}
\put(20,30){\line(1,-1){10}}
\put(15,-3){\makebox(0,0)[]{Kh15: (8,7,7,6,6,6,5,5,4)}}
\end{picture}
}

\end{picture}

\begin{picture}(110,40)(0,0)

\put(0,0){
\begin{picture}(30,40)(0,-5)
\thicklines
\multiput(0,0)(0,10){4}{\line(1,0){30}}
\multiput(0,0)(10,0){4}{\line(0,1){30}}
\multiput(0,0)(10,0){4}{\circle*{1}}
\multiput(0,10)(10,0){4}{\circle*{1}}
\multiput(0,20)(10,0){4}{\circle*{1}}
\multiput(0,30)(10,0){4}{\circle*{1}}
\put(0,0){\line(1,1){10}}
\put(10,10){\line(1,-1){10}}
\put(20,0){\line(1,1){10}}
\put(0,20){\line(1,-1){10}}
\put(10,10){\line(1,1){10}}
\put(20,20){\line(1,-1){10}}
\put(0,20){\line(1,1){10}}
\put(10,30){\line(1,-1){10}}
\put(20,20){\line(1,1){10}}
\put(15,-3){\makebox(0,0)[]{Kh16: (8,8,8,6,6,6,4,4,4)}}
\end{picture}
}

\put(40,0){
\begin{picture}(30,40)(0,-5)
\thicklines
\multiput(0,0)(0,10){4}{\line(1,0){30}}
\multiput(0,0)(10,0){4}{\line(0,1){30}}
\multiput(0,0)(10,0){4}{\circle*{1}}
\multiput(0,10)(10,0){4}{\circle*{1}}
\multiput(0,20)(10,0){4}{\circle*{1}}
\multiput(0,30)(10,0){4}{\circle*{1}}
\put(0,10){\line(1,-1){10}}
\put(10,10){\line(1,-1){10}}
\put(20,0){\line(1,1){10}}
\put(0,20){\line(1,-1){10}}
\put(10,10){\line(1,1){10}}
\put(20,20){\line(1,-1){10}}
\put(0,20){\line(1,1){10}}
\put(10,30){\line(1,-1){10}}
\put(20,30){\line(1,-1){10}}
\put(15,-3){\makebox(0,0)[]{Kh17: (8,7,7,7,7,6,4,4,4)}}
\end{picture}
}

\put(80,0){
\begin{picture}(30,40)(0,-5)
\thicklines
\multiput(0,0)(0,10){4}{\line(1,0){30}}
\multiput(0,0)(10,0){4}{\line(0,1){30}}
\multiput(0,0)(10,0){4}{\circle*{1}}
\multiput(0,10)(10,0){4}{\circle*{1}}
\multiput(0,20)(10,0){4}{\circle*{1}}
\multiput(0,30)(10,0){4}{\circle*{1}}
\put(0,0){\line(1,1){10}}
\put(10,10){\line(1,-1){10}}
\put(20,0){\line(1,1){10}}
\put(0,20){\line(1,-1){10}}
\put(10,10){\line(1,1){10}}
\put(20,20){\line(1,-1){10}}
\put(0,20){\line(1,1){10}}
\put(10,30){\line(1,-1){10}}
\put(20,30){\line(1,-1){10}}
\put(15,-3){\makebox(0,0)[]{Kh18: (8,8,7,7,6,5,5,4,4)}}
\end{picture}
}

\end{picture}

\begin{picture}(110,40)(0,0)

\put(0,0){
\begin{picture}(30,40)(0,-5)
\thicklines
\multiput(0,0)(0,10){4}{\line(1,0){30}}
\multiput(0,0)(10,0){4}{\line(0,1){30}}
\multiput(0,0)(10,0){4}{\circle*{1}}
\multiput(0,10)(10,0){4}{\circle*{1}}
\multiput(0,20)(10,0){4}{\circle*{1}}
\multiput(0,30)(10,0){4}{\circle*{1}}
\put(0,10){\line(1,-1){10}}
\put(10,0){\line(1,1){10}}
\put(20,10){\line(1,-1){10}}
\put(0,20){\line(1,-1){10}}
\put(10,10){\line(1,1){10}}
\put(20,20){\line(1,-1){10}}
\put(0,20){\line(1,1){10}}
\put(10,30){\line(1,-1){10}}
\put(20,20){\line(1,1){10}}
\put(15,-3){\makebox(0,0)[]{Kh19: (8,8,7,6,6,6,5,4,4)}}
\end{picture}
}

\put(40,0){
\begin{picture}(30,40)(0,-5)
\thicklines
\multiput(0,0)(0,10){4}{\line(1,0){30}}
\multiput(0,0)(10,0){4}{\line(0,1){30}}
\multiput(0,0)(10,0){4}{\circle*{1}}
\multiput(0,10)(10,0){4}{\circle*{1}}
\multiput(0,20)(10,0){4}{\circle*{1}}
\multiput(0,30)(10,0){4}{\circle*{1}}
\put(0,10){\line(1,-1){10}}
\put(10,0){\line(1,1){10}}
\put(20,0){\line(1,1){10}}
\put(0,20){\line(1,-1){10}}
\put(10,10){\line(1,1){10}}
\put(20,20){\line(1,-1){10}}
\put(0,20){\line(1,1){10}}
\put(10,30){\line(1,-1){10}}
\put(20,20){\line(1,1){10}}
\put(15,-3){\makebox(0,0)[]{Kh20: (8,8,8,6,5,5,5,5,4)}}
\end{picture}
}

\put(80,0){
\begin{picture}(30,40)(0,-5)
\thicklines
\multiput(0,0)(0,10){4}{\line(1,0){30}}
\multiput(0,0)(10,0){4}{\line(0,1){30}}
\multiput(0,0)(10,0){4}{\circle*{1}}
\multiput(0,10)(10,0){4}{\circle*{1}}
\multiput(0,20)(10,0){4}{\circle*{1}}
\multiput(0,30)(10,0){4}{\circle*{1}}
\put(0,10){\line(1,-1){10}}
\put(10,10){\line(1,-1){10}}
\put(20,10){\line(1,-1){10}}
\put(0,20){\line(1,-1){10}}
\put(10,10){\line(1,1){10}}
\put(20,20){\line(1,-1){10}}
\put(0,20){\line(1,1){10}}
\put(10,30){\line(1,-1){10}}
\put(20,20){\line(1,1){10}}
\put(15,-3){\makebox(0,0)[]{Kh21: (8,7,7,7,6,5,5,5,4)}}
\end{picture}
}

\end{picture}
\caption{Irreducible triangulations of the Klein bottle, Kh13 - Kh21}
\label{fig2}
\end{center}
\end{figure}
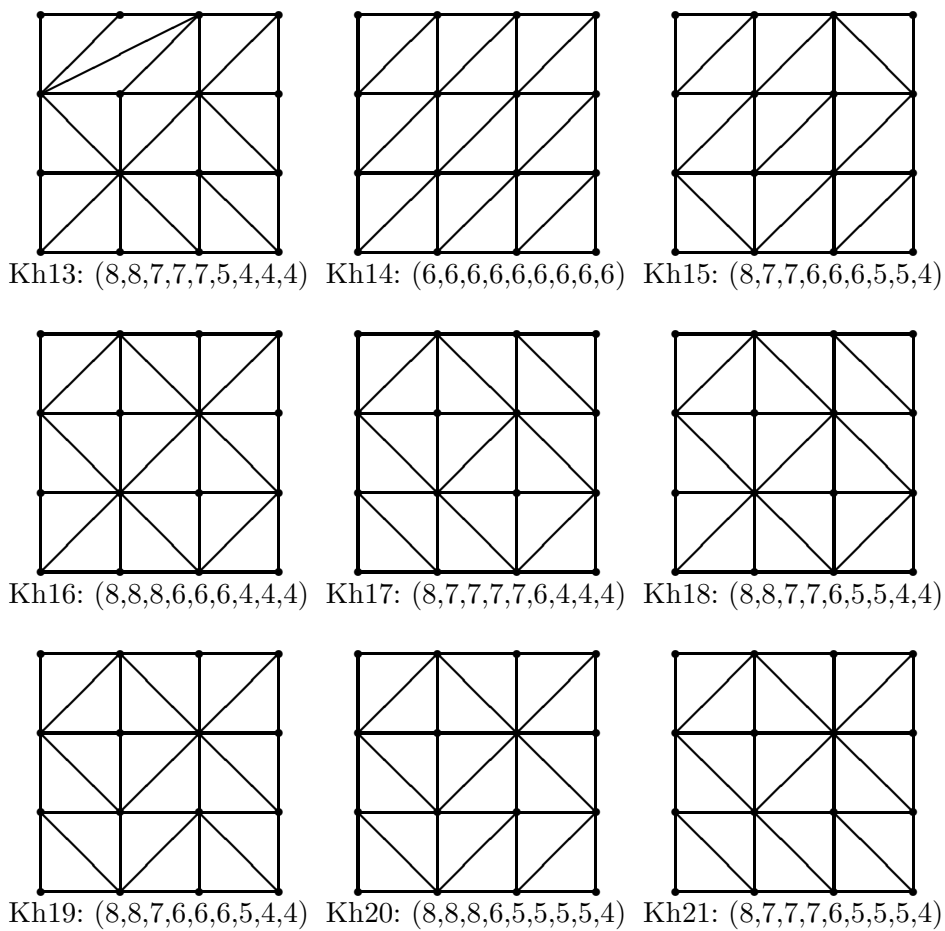

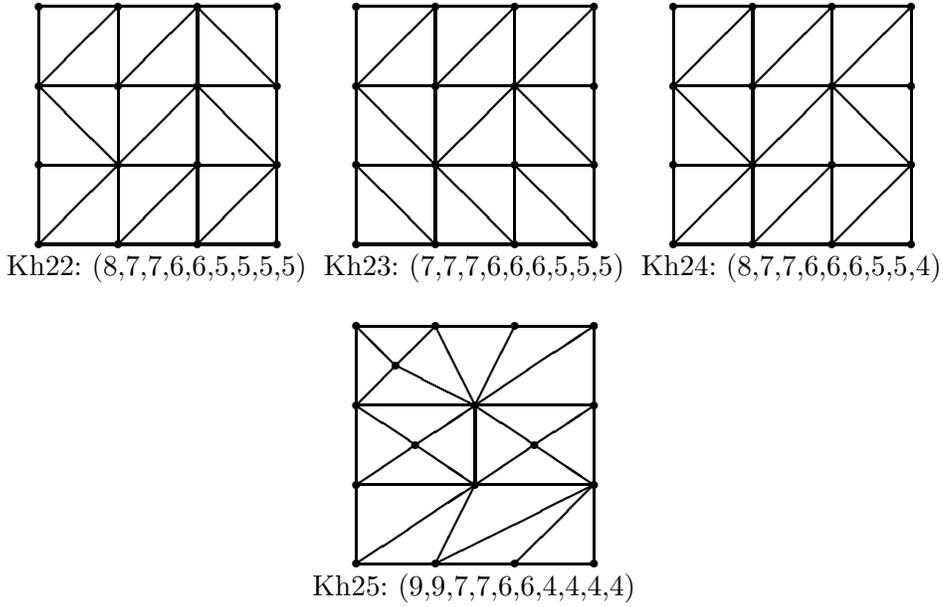
\begin{figure}
\begin{center}
\setlength{\unitlength}{3pt}
\small

\begin{picture}(110,40)(0,0)

\put(0,0){
\begin{picture}(30,40)(0,-5)
\thicklines
\multiput(0,0)(0,10){4}{\line(1,0){30}}
\multiput(0,0)(10,0){4}{\line(0,1){30}}
\multiput(0,0)(10,0){4}{\circle*{1}}
\multiput(0,10)(10,0){4}{\circle*{1}}
\multiput(0,20)(10,0){4}{\circle*{1}}
\multiput(0,30)(10,0){4}{\circle*{1}}
\put(0,0){\line(1,1){10}}
\put(10,0){\line(1,1){10}}
\put(20,0){\line(1,1){10}}
\put(0,20){\line(1,-1){10}}
\put(10,10){\line(1,1){10}}
\put(20,20){\line(1,-1){10}}
\put(0,20){\line(1,1){10}}
\put(10,20){\line(1,1){10}}
\put(20,30){\line(1,-1){10}}
\put(15,-3){\makebox(0,0)[]{Kh22: (8,7,7,6,6,5,5,5,5)}}
\end{picture}
}

\put(40,0){
\begin{picture}(30,40)(0,-5)
\thicklines
\multiput(0,0)(0,10){4}{\line(1,0){30}}
\multiput(0,0)(10,0){4}{\line(0,1){30}}
\multiput(0,0)(10,0){4}{\circle*{1}}
\multiput(0,10)(10,0){4}{\circle*{1}}
\multiput(0,20)(10,0){4}{\circle*{1}}
\multiput(0,30)(10,0){4}{\circle*{1}}
\put(0,10){\line(1,-1){10}}
\put(10,10){\line(1,-1){10}}
\put(20,10){\line(1,-1){10}}
\put(0,20){\line(1,-1){10}}
\put(10,10){\line(1,1){10}}
\put(20,20){\line(1,-1){10}}
\put(0,20){\line(1,1){10}}
\put(10,20){\line(1,1){10}}
\put(20,20){\line(1,1){10}}
\put(15,-3){\makebox(0,0)[]{Kh23: (7,7,7,6,6,6,5,5,5)}}
\end{picture}
}

\put(80,0){
\begin{picture}(30,40)(0,-5)
\thicklines
\multiput(0,0)(0,10){4}{\line(1,0){30}}
\multiput(0,0)(10,0){4}{\line(0,1){30}}
\multiput(0,0)(10,0){4}{\circle*{1}}
\multiput(0,10)(10,0){4}{\circle*{1}}
\multiput(0,20)(10,0){4}{\circle*{1}}
\multiput(0,30)(10,0){4}{\circle*{1}}
\put(0,0){\line(1,1){10}}
\put(10,0){\line(1,1){10}}
\put(20,0){\line(1,1){10}}
\put(0,20){\line(1,-1){10}}
\put(10,10){\line(1,1){10}}
\put(20,20){\line(1,-1){10}}
\put(0,20){\line(1,1){10}}
\put(10,20){\line(1,1){10}}
\put(20,20){\line(1,1){10}}
\put(15,-3){\makebox(0,0)[]{Kh24: (8,7,7,6,6,6,5,5,4)}}
\end{picture}
}

\end{picture}

\begin{picture}(110,40)(0,0)

\put(40,0){
\begin{picture}(30,40)(0,-5)
\thicklines
\multiput(0,0)(0,10){4}{\line(1,0){30}}
\multiput(0,0)(30,0){2}{\line(0,1){30}}
\multiput(0,0)(10,0){4}{\circle*{1}}
\multiput(0,10)(15,0){3}{\circle*{1}}
\multiput(7.5,15)(15,0){2}{\circle*{1}}
\multiput(0,20)(15,0){3}{\circle*{1}}
\put(5,25){\circle*{1}}
\multiput(0,30)(10,0){4}{\circle*{1}}
\put(0,0){\line(3,2){15}}
\put(10,0){\line(1,2){5}}
\put(10,0){\line(2,1){20}}
\put(20,0){\line(1,1){10}}
\put(0,10){\line(3,2){15}}
\put(0,20){\line(3,-2){15}}
\put(15,10){\line(0,1){10}}
\put(15,10){\line(3,2){15}}
\put(15,20){\line(3,-2){15}}
\put(0,20){\line(1,1){10}}
\put(0,30){\line(1,-1){5}}
\put(5,25){\line(2,-1){10}}
\put(10,30){\line(1,-2){5}}
\put(15,20){\line(1,2){5}}
\put(15,20){\line(3,2){15}}
\put(15,-3){\makebox(0,0)[]{Kh25: (9,9,7,7,6,6,4,4,4,4)}}
\end{picture}
}

\end{picture}
\caption{Irreducible triangulations of the Klein bottle, Kh22 - Kh25}
\label{fig3}
\end{center}
\end{figure}

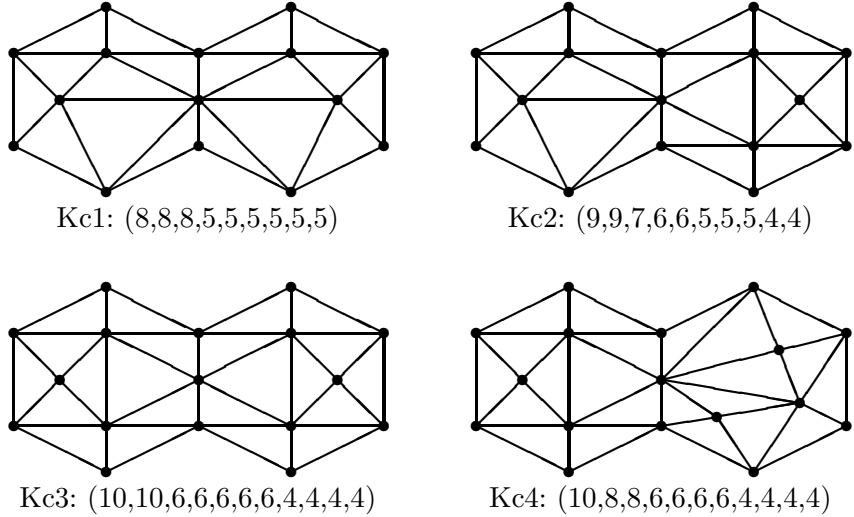
\begin{figure}
\begin{center}
\setlength{\unitlength}{3.5pt}
\small

\begin{picture}(90,30)(0,0)

\put(0,0){
\begin{picture}(40,25)(-20,-15)
\thicklines
\multiput(-20,-5)(20,0){3}{\line(0,1){10}}
\put(-20,-5){\line(2,-1){10}}
\put(-10,-10){\line(2,1){10}}
\put(0,-5){\line(2,-1){10}}
\put(10,-10){\line(2,1){10}}
\put(-20,5){\line(2,1){10}}
\put(-10,10){\line(2,-1){10}}
\put(0,5){\line(2,1){10}}
\put(10,10){\line(2,-1){10}}
\multiput(-10,-10)(20,0){2}{\circle*{1}}
\multiput(-20,-5)(20,0){3}{\circle*{1}}
\put(0,0){\circle*{1}}
\multiput(-20,5)(20,0){3}{\circle*{1}}
\multiput(-10,10)(20,0){2}{\circle*{1}}

\put(-15,0){\circle*{1}}
\put(-10,5){\circle*{1}}
\put(-10,-10){\line(-1,2){5}}
\put(-10,-10){\line(1,1){10}}
\put(-20,-5){\line(1,1){5}}
\put(-15,0){\line(1,0){15}}
\put(-15,0){\line(-1,1){5}}
\put(-15,0){\line(1,1){5}}
\put(-10,5){\line(2,-1){10}}
\put(-20,5){\line(1,0){20}}
\put(-10,5){\line(0,1){5}}

\put(15,0){\circle*{1}}
\put(10,5){\circle*{1}}
\put(10,-10){\line(1,2){5}}
\put(10,-10){\line(-1,1){10}}
\put(20,-5){\line(-1,1){5}}
\put(15,0){\line(-1,0){15}}
\put(15,0){\line(1,1){5}}
\put(15,0){\line(-1,1){5}}
\put(10,5){\line(-2,-1){10}}
\put(20,5){\line(-1,0){20}}
\put(10,5){\line(0,1){5}}
\put(0,-13){\makebox(0,0)[]{Kc1: (8,8,8,5,5,5,5,5,5)}}
\end{picture}
}

\put(50,0){
\begin{picture}(40,25)(-20,-15)
\thicklines
\multiput(-20,-5)(20,0){3}{\line(0,1){10}}
\put(-20,-5){\line(2,-1){10}}
\put(-10,-10){\line(2,1){10}}
\put(0,-5){\line(2,-1){10}}
\put(10,-10){\line(2,1){10}}
\put(-20,5){\line(2,1){10}}
\put(-10,10){\line(2,-1){10}}
\put(0,5){\line(2,1){10}}
\put(10,10){\line(2,-1){10}}
\multiput(-10,-10)(20,0){2}{\circle*{1}}
\multiput(-20,-5)(20,0){3}{\circle*{1}}
\put(0,0){\circle*{1}}
\multiput(-20,5)(20,0){3}{\circle*{1}}
\multiput(-10,10)(20,0){2}{\circle*{1}}

\put(-15,0){\circle*{1}}
\put(-10,5){\circle*{1}}
\put(-10,-10){\line(-1,2){5}}
\put(-10,-10){\line(1,1){10}}
\put(-20,-5){\line(1,1){5}}
\put(-15,0){\line(1,0){15}}
\put(-15,0){\line(-1,1){5}}
\put(-15,0){\line(1,1){5}}
\put(-10,5){\line(2,-1){10}}
\put(-20,5){\line(1,0){20}}
\put(-10,5){\line(0,1){5}}

\put(10,-5){\circle*{1}}
\put(15,0){\circle*{1}}
\put(10,5){\circle*{1}}
\put(10,-10){\line(0,1){20}}
\put(0,-5){\line(1,0){20}}
\put(10,-5){\line(-2,1){10}}
\put(10,-5){\line(1,1){10}}
\put(0,0){\line(2,1){10}}
\put(10,5){\line(1,-1){10}}
\put(20,5){\line(-1,0){20}}
\put(0,-13){\makebox(0,0)[]{Kc2: (9,9,7,6,6,5,5,5,4,4)}}
\end{picture}
}

\end{picture}

\begin{picture}(90,30)(0,0)

\put(0,0){
\begin{picture}(40,25)(-20,-15)
\thicklines
\multiput(-20,-5)(20,0){3}{\line(0,1){10}}
\put(-20,-5){\line(2,-1){10}}
\put(-10,-10){\line(2,1){10}}
\put(0,-5){\line(2,-1){10}}
\put(10,-10){\line(2,1){10}}
\put(-20,5){\line(2,1){10}}
\put(-10,10){\line(2,-1){10}}
\put(0,5){\line(2,1){10}}
\put(10,10){\line(2,-1){10}}
\multiput(-10,-10)(20,0){2}{\circle*{1}}
\multiput(-20,-5)(20,0){3}{\circle*{1}}
\put(0,0){\circle*{1}}
\multiput(-20,5)(20,0){3}{\circle*{1}}
\multiput(-10,10)(20,0){2}{\circle*{1}}

\put(-10,-5){\circle*{1}}
\put(-15,0){\circle*{1}}
\put(-10,5){\circle*{1}}
\put(-10,-10){\line(0,1){20}}
\put(0,-5){\line(-1,0){20}}
\put(-10,-5){\line(2,1){10}}
\put(-10,-5){\line(-1,1){10}}
\put(0,0){\line(-2,1){10}}
\put(-10,5){\line(-1,-1){10}}
\put(-20,5){\line(1,0){20}}

\put(10,-5){\circle*{1}}
\put(15,0){\circle*{1}}
\put(10,5){\circle*{1}}
\put(10,-10){\line(0,1){20}}
\put(0,-5){\line(1,0){20}}
\put(10,-5){\line(-2,1){10}}
\put(10,-5){\line(1,1){10}}
\put(0,0){\line(2,1){10}}
\put(10,5){\line(1,-1){10}}
\put(20,5){\line(-1,0){20}}
\put(0,-13){\makebox(0,0)[]{Kc3: (10,10,6,6,6,6,6,4,4,4,4)}}
\end{picture}
}

\put(50,0){
\begin{picture}(40,25)(-20,-15)
\thicklines
\multiput(-20,-5)(20,0){3}{\line(0,1){10}}
\put(-20,-5){\line(2,-1){10}}
\put(-10,-10){\line(2,1){10}}
\put(0,-5){\line(2,-1){10}}
\put(10,-10){\line(2,1){10}}
\put(-20,5){\line(2,1){10}}
\put(-10,10){\line(2,-1){10}}
\put(0,5){\line(2,1){10}}
\put(10,10){\line(2,-1){10}}
\multiput(-10,-10)(20,0){2}{\circle*{1}}
\multiput(-20,-5)(20,0){3}{\circle*{1}}
\put(0,0){\circle*{1}}
\multiput(-20,5)(20,0){3}{\circle*{1}}
\multiput(-10,10)(20,0){2}{\circle*{1}}

\put(-10,-5){\circle*{1}}
\put(-15,0){\circle*{1}}
\put(-10,5){\circle*{1}}
\put(-10,-10){\line(0,1){20}}
\put(0,-5){\line(-1,0){20}}
\put(-10,-5){\line(2,1){10}}
\put(-10,-5){\line(-1,1){10}}
\put(0,0){\line(-2,1){10}}
\put(-10,5){\line(-1,-1){10}}
\put(-20,5){\line(1,0){20}}

\put(6,-4){\circle*{1}}
\put(15,-2.5){\circle*{1}}
\put(12.73,3.18){\circle*{1}}
\put(10,-10){\line(-2,3){4}}
\put(10,-10){\line(2,3){10}}
\put(0,-5){\line(6,1){15}}
\put(20,-5){\line(-2,1){5}}
\put(0,0){\line(3,-2){6}}
\put(0,0){\line(6,-1){15}}
\put(0,0){\line(4,1){20}}
\put(0,0){\line(1,1){10}}
\put(10,10){\line(2,-5){5}}
\put(0,-13){\makebox(0,0)[]{Kc4: (10,8,8,6,6,6,6,4,4,4,4)}}
\end{picture}
}

\end{picture}

\caption{Irreducible triangulations of the Klein bottle, Kc1 - Kc4}
\label{fig4}
\end{center}
\end{figure}

To show that there are at least 29 irreducible triangulations of the 
Klein bottle we exhibit them.
Figures~\ref{fig1},~\ref{fig2},~\ref{fig3}, and~\ref{fig4} 
show the complete list.
Figures~\ref{fig1},~\ref{fig2}, and~\ref{fig4} are the same as 
Figs. 13, 14, and 15 of~\cite{LN}.
In Figs.~\ref{fig1},~\ref{fig2}, and~\ref{fig3} the pair of horizontal 
sides of each rectangle are identified in parallel and the pair of
vertical sides are identified in
antiparallel to obtain an actual triangulation of the Klein bottle.
In~\cite{LN} such triangulations are classified as 
{\em handle type}.
In Fig.~\ref{fig4} each triangulation 
can be obtained from two copies of irreducible triangulations of the
projective plane by removing one triangular face from each copy and 
pasting them together along boundaries of the removed faces.
For each of the two hexagons in each graph
identify each antipodal pair of vertices on the boundary to obtain the 
triangulation of the projective plane less one face.  
When these identifications have been made for both parts of a graph a 
triangulation of the Klein bottle is produced.
This type of triangulation is classified as {\em crosscap type} in~\cite{LN}.

The irreducible triangulations shown in Fig.~\ref{fig3}, which are denoted
Kh22, Kh23, Kh24, and Kh25, were not listed in \cite{LN}.  
It can be seen that 
they are indeed irreducible by checking that each edge $e$ 
is contained in a cycle of 
length 3 other than the two  which bound the faces incident to $e$. 

The triangulations Kh22 through Kh25 are not equivalent to any of the other
triangulations.  
In fact, they are not isomorphic as graphs to any of the other 
triangulations.
The degree sequence of each triangulation is shown in Figs.~\ref{fig1} 
through~\ref{fig4}.  
The degree sequences of Kh23 and Kh25 are unique.  
The degree sequence of Kh22 and Kh10 are the same. However, in Kh10
the two vertices of degree 6 are adjacent, while in Kh22 
the two vertices of degree 6 are not adjacent.  
Similarly, the degree sequences of Kh24 and Kh15 are the same. In Kh15
the two vertices of degree 5 are adjacent, while in Kh24 
the two vertices of degree 5 are not adjacent.

\section{Construction of Kh25}
\label{mod-Kh25}

The construction of Lawrencenko and Negami can be modified to produce
all~29 
irreducible triangulations of the Klein bottle and to show that there 
are no others.
We do not need to modify the basic proof of Lemma 6 of~\cite{LN}.
Instead we correct some of the details of their proof.
Assuming that the reader is familiar with the proof in \cite{LN} 
we do not give its details here.

The triangulation Kh25 was missed in determining the partial structures 
inside the rectangle $R_{01}$ as defined in \cite{LN}.
In {\em Step 2. Recognizing the inside of $R_{01}$} of~\cite{LN},
in the last full paragraph on Page 279, Line -11 the authors incorrectly 
stated that
``there is no chord incident to $y_1$ outside the polygon''
$W_2 \cup y_1 b' a' x_1 x_2$
and used this conclusion to the end of that paragraph.
However, one needs to consider both the presence and the absence of 
such a chord.

First assume that vertex $y_1$ is not adjacent to $a'$ or to $a''$.
Then there is no chord incident to $y_1$ outside the polygon 
$W_2 \cup y_1 b' a' x_1 x_2$.
The path $W_1$ has length $1$, that is, $W_1 = x_1 y_1$ 
by Lemma 3 of \cite{LN}.
Similarly, $W_2 = x_2 y_1$ and each of the two quadrilateral 
regions bounded by 
$a' b' y_1 x_1$ and $a'' c'' y_1 x_2$ contains only one diagonal by 
Lemma 3 of \cite{LN}.
In this case, by assumption $y_1 a'$ and $y_1 a''$ are not diagonals, 
thus we have the partial structure $R_{01}\!-\!6$.

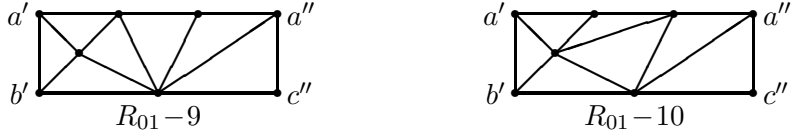
\begin{figure}
\begin{center}
\setlength{\unitlength}{3pt}
\small

\begin{picture}(100,16)(0,0)
\put(0,0){

\begin{picture}(40,16)(0,-5)
\thicklines
\multiput(0,0)(0,10){2}{\line(1,0){30}}
\multiput(0,0)(30,0){2}{\line(0,1){10}}
\multiput(0,0)(15,0){3}{\circle*{1}}
\multiput(0,10)(10,0){4}{\circle*{1}}
\put(5,5){\circle*{1}}
\put(0,0){\line(1,1){10}}
\put(0,10){\line(1,-1){5}}
\put(5,5){\line(2,-1){10}}
\put(10,10){\line(1,-2){5}}
\put(15,0){\line(1,2){5}}
\put(15,0){\line(3,2){15}}
\put(0,10){\makebox(0,0)[r]{$a'$\ }}
\put(0,0){\makebox(0,0)[r]{$b'$\ }}
\put(30,10){\makebox(0,0)[l]{\ $a''$}}
\put(30,0){\makebox(0,0)[l]{\ $c''$}}

\put(15,-3){\makebox(0,0)[]{$R_{01}\!-\!9$}}
\end{picture}
}

\put(60,0){

\begin{picture}(40,16)(0,-5)
\thicklines
\multiput(0,0)(0,10){2}{\line(1,0){30}}
\multiput(0,0)(30,0){2}{\line(0,1){10}}
\multiput(0,0)(15,0){3}{\circle*{1}}
\multiput(0,10)(10,0){4}{\circle*{1}}
\put(5,5){\circle*{1}}
\put(0,0){\line(1,1){10}}
\put(0,10){\line(1,-1){5}}
\put(5,5){\line(2,-1){10}}
\put(20,10){\line(-3,-1){15}}
\put(15,0){\line(1,2){5}}
\put(15,0){\line(3,2){15}}
\put(0,10){\makebox(0,0)[r]{$a'$\ }}
\put(0,0){\makebox(0,0)[r]{$b'$\ }}
\put(30,10){\makebox(0,0)[l]{\ $a''$}}
\put(30,0){\makebox(0,0)[l]{\ $c''$}}

\put(15,-3){\makebox(0,0)[]{$R_{01}\!-\!10$}}
\end{picture}
}

\end{picture}

\caption{Additional structures inside $R_{01}$.}
\label{R9}
\end{center}
\end{figure}

Now by symmetry, without loss of generality, assume that $y_1$ is 
adjacent to $a''$.
The edge $x_2 y_1$ is not a chord outside the polygon 
$W_1 \cup x_1 x_2 a'' y_1$, thus $W_2$ has length $1$, that is, 
$W_2 = x_2 y_1$.
If $W_1$ has length $1$, 
then the polygon $y_1 b' a' x_1$ either has a diagonal 
$b' x_1$ and we have the partial structure $R_{01}\!-\!5$, or 
there is one internal
vertex in $y_1 b' a' x_1$ and we have an additional partial structure 
$R_{01}\!-\!9$ which is shown in Fig.~\ref{R9}.
The polygon $y_1 b' a' x_1$ in $R_{01}\!-\!9$ has an internal vertex and 
$x_1 b''$ must be an external chord in $R_{02}$.
If, on the other hand, $W_1$ has length greater than $1$, then there is exactly
one vertex $z$ on $W_1$ and inside the polygon $y_1 b' a' x_1 x_2$ 
by Lemma 3 of \cite{LN}.
In this case we show that $z$ is adjacent to all the vertices of the polygon 
$y_1 b' a' x_1 x_2$.
Since $W_1 = x_1 z y_1$ is a minimum length path joining $x_1$ and $y_1$, the
edge $x_1 y_1$ is not a diagonal of the polygon $x_1 z y_1 x_2$.
There is no interior vertex in the polygon $x_1 z y_1 x_2$, hence $z x_2$ is a
diagonal.
The polygon $y_1 b' a' x_1 z$ has no interior vertex and neither $y_1 a'$ nor 
$y_1 x_1$ are edges, hence $b' z$ is an edge in $y_1 b' a' x_1 z$.
If $b' x_1$ is an edge in $R_{01}$, then the polygon $b' x_1 x_2 y_1$ would 
contain vertex $z$; 
hence by Lemma 3 of \cite{LN}, there would be a chord $x_1 y_1$
outside the polygon $b' x_1 x_2 y_1$, which is impossible.
Therefore, $b' x_1$ is not an edge in $R_{01}$ and $z a'$ is the diagonal of 
the polygon $b' a' x_1 z$.
In this case we have another additional partial structure $R_{01}\!-\!10$ 
which is shown in Fig.~\ref{R9}.
The polygon $y_1 b' a' x_1 x_2$ in $R_{01}\!-\!10$ has an internal vertex, and 
$x_1 b''$ must be an external chord in $R_{02}$.

In {\em Step 4. Composing partial structures in triangulations\/} 
of~\cite{LN}, the additional configurations using the
partial structures $R_{01}\!-\!9$ and $R_{01}\!-\!10$ must be classified.
Both $R_{01}\!-\!9$ and $R_{01}\!-\!10$ require 
the edge $x_1 b''$ in $R_{02}$,
thus the additional configurations are
$\{[i,1,-k]:i=9,10;k=1,2,3\}$ and 
$\{[i,j,-k]:i=9,10;j=2,3;k=7,8\}$.
The configuration $[9,3,-8]$ is added to group (i) and is the additional 
irreducible triangulation Kh25.
The configuration $[9,3,-7]$ is equivalent to $[9,3,-8]$.
The configuration $[10,1,-1]$, which is equivalent to $[5,2,-5]$,
is contrary to assumption (III).
The configurations $[9,1,-1]$,
$[9,1,-2]$, and
$[9,1,-3]$ are added to group (v).
The configurations $[9,2,-7]$, $[9,2,-8]$,
$[10,1,-2]$, $[10,1,-3]$,
$[10,2,-7]$, $[10,2,-8]$,
$[10,3,-7]$, and $[10,3,-8]$
are added to group (vii).

\section{Construction of Kh22, Kh23, and Kh24}
\label{mod-Kh222324}

In {\em Step 5.  Classifying triangulations up to equivalence\/} 
(the first paragraph on page 283) of the proof of Lemma 6 of~\cite{LN},
diagonals were added to the partial structures and the results were 
classified up to equivalence.  
No details were given regarding how this classification was accomplished.
Only 8 irreducible triangulations (Kh14 through Kh21) which can be
obtained by adding diagonals to the partial structure PS1 (Fig. 8 of~\cite{LN})
were listed in~\cite{LN}.
There are 3 additional irreducible triangulations (Kh22 through Kh24) which
can be obtained from PS1.
We examine in detail how to add diagonals to PS1 and classify all the resulting
triangulations.
This procedure requires that we check only $16$, not all
$2^9$, configurations obtained by adding diagonals to the nine quadrilateral 
regions of PS1.

\begin{figure}
\begin{center}
\setlength{\unitlength}{1.5pt}
\small

\begin{picture}(230,40)(0,0)

\put(0,0){
\begin{picture}(30,40)(0,-5)
\thicklines
\multiput(0,0)(0,10){4}{\line(1,0){30}}
\multiput(0,0)(10,0){4}{\line(0,1){30}}
\multiput(0,0)(10,0){4}{\circle*{1}}
\multiput(0,10)(10,0){4}{\circle*{1}}
\multiput(0,20)(10,0){4}{\circle*{1}}
\multiput(0,30)(10,0){4}{\circle*{1}}
\put(0,10){\line(1,1){10}}
\put(10,10){\line(1,1){10}}
\put(20,20){\line(1,-1){10}}
\end{picture}
}

\put(40,0){
\begin{picture}(30,40)(0,-5)
\thicklines
\multiput(0,0)(0,10){4}{\line(1,0){30}}
\multiput(0,0)(10,0){4}{\line(0,1){30}}
\multiput(0,0)(10,0){4}{\circle*{1}}
\multiput(0,10)(10,0){4}{\circle*{1}}
\multiput(0,20)(10,0){4}{\circle*{1}}
\multiput(0,30)(10,0){4}{\circle*{1}}
\put(0,10){\line(1,1){10}}
\put(10,20){\line(1,-1){10}}
\put(20,20){\line(1,-1){10}}
\end{picture}
}

\put(80,0){
\begin{picture}(30,40)(0,-5)
\thicklines
\multiput(0,0)(0,10){4}{\line(1,0){30}}
\multiput(0,0)(10,0){4}{\line(0,1){30}}
\multiput(0,0)(10,0){4}{\circle*{1}}
\multiput(0,10)(10,0){4}{\circle*{1}}
\multiput(0,20)(10,0){4}{\circle*{1}}
\multiput(0,30)(10,0){4}{\circle*{1}}
\put(0,20){\line(1,-1){10}}
\put(10,20){\line(1,-1){10}}
\put(20,20){\line(1,-1){10}}
\end{picture}
}

\put(120,0){
\begin{picture}(30,40)(0,-5)
\thicklines
\multiput(0,0)(0,10){4}{\line(1,0){30}}
\multiput(0,0)(10,0){4}{\line(0,1){30}}
\multiput(0,0)(10,0){4}{\circle*{1}}
\multiput(0,10)(10,0){4}{\circle*{1}}
\multiput(0,20)(10,0){4}{\circle*{1}}
\multiput(0,30)(10,0){4}{\circle*{1}}
\put(0,20){\line(1,-1){10}}
\put(10,20){\line(1,-1){10}}
\put(20,10){\line(1,1){10}}
\end{picture}
}

\put(160,0){
\begin{picture}(30,40)(0,-5)
\thicklines
\multiput(0,0)(0,10){4}{\line(1,0){30}}
\multiput(0,0)(10,0){4}{\line(0,1){30}}
\multiput(0,0)(10,0){4}{\circle*{1}}
\multiput(0,10)(10,0){4}{\circle*{1}}
\multiput(0,20)(10,0){4}{\circle*{1}}
\multiput(0,30)(10,0){4}{\circle*{1}}
\put(0,20){\line(1,-1){10}}
\put(10,10){\line(1,1){10}}
\put(20,10){\line(1,1){10}}
\end{picture}
}

\put(200,0){
\begin{picture}(30,40)(0,-5)
\thicklines
\multiput(0,0)(0,10){4}{\line(1,0){30}}
\multiput(0,0)(10,0){4}{\line(0,1){30}}
\multiput(0,0)(10,0){4}{\circle*{1}}
\multiput(0,10)(10,0){4}{\circle*{1}}
\multiput(0,20)(10,0){4}{\circle*{1}}
\multiput(0,30)(10,0){4}{\circle*{1}}
\put(0,10){\line(1,1){10}}
\put(10,10){\line(1,1){10}}
\put(20,10){\line(1,1){10}}
\put(15,-5){\makebox(0,0)[]{PS1.1}}
\end{picture}
}

\end{picture}

\begin{picture}(230,45)(0,0)

\put(160,0){
\begin{picture}(30,40)(0,-5)
\thicklines
\multiput(0,0)(0,10){4}{\line(1,0){30}}
\multiput(0,0)(10,0){4}{\line(0,1){30}}
\multiput(0,0)(10,0){4}{\circle*{1}}
\multiput(0,10)(10,0){4}{\circle*{1}}
\multiput(0,20)(10,0){4}{\circle*{1}}
\multiput(0,30)(10,0){4}{\circle*{1}}
\put(0,10){\line(1,1){10}}
\put(10,20){\line(1,-1){10}}
\put(20,10){\line(1,1){10}}
\end{picture}
}

\put(200,0){
\begin{picture}(30,40)(0,-5)
\thicklines
\multiput(0,0)(0,10){4}{\line(1,0){30}}
\multiput(0,0)(10,0){4}{\line(0,1){30}}
\multiput(0,0)(10,0){4}{\circle*{1}}
\multiput(0,10)(10,0){4}{\circle*{1}}
\multiput(0,20)(10,0){4}{\circle*{1}}
\multiput(0,30)(10,0){4}{\circle*{1}}
\put(0,20){\line(1,-1){10}}
\put(10,10){\line(1,1){10}}
\put(20,20){\line(1,-1){10}}
\put(15,-5){\makebox(0,0)[]{PS1.2}}
\end{picture}
}

\end{picture}

\caption{PS1 with diagonals in $R_{12}$.}
\label{middle row}
\end{center}
\end{figure}
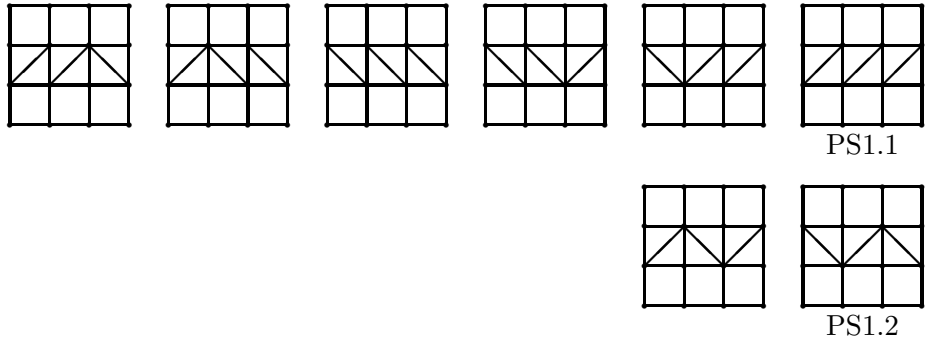

Consider the possible configurations for diagonals 
in $R_{12}$ of PS1.  
Figure~\ref{middle row} shows all eight partial structures obtained 
by adding diagonals to $R_{12}$ of PS1.
Going from left to right, each partial structure is transformed into the
next by removing the left column, reflecting it vertically,
and pasting it onto the right side.  
Thus we can assume that PS1 has one of the two rightmost
partial structures, PS1.1 or PS1.2, of Fig.~\ref{middle row}.

\begin{figure}
\begin{center}
\setlength{\unitlength}{3pt}

\begin{picture}(90,45)(0,0)


\put(0,0){
\begin{picture}(40,45)(-5,-10)
\thicklines
\multiput(0,0)(10,0){4}{\circle*{1}}
\multiput(0,10)(10,0){4}{\circle*{1}}
\multiput(0,20)(10,0){4}{\circle*{1}}
\multiput(0,30)(10,0){4}{\circle*{1}}
\multiput(0,0)(0,10){4}{\line(1,0){30}}
\multiput(0,0)(10,0){4}{\line(0,1){30}}
\put(00,10){\line(1,1){10}}
\put(10,10){\line(1,1){10}}
\put(20,10){\line(1,1){10}}
\put(0,30){\makebox(0,0)[r]{$a$\ }}
\put(0,20){\makebox(0,0)[r]{$b$\ }}
\put(0,10){\makebox(0,0)[r]{$c$\ }}
\put(0,00){\makebox(0,0)[r]{$a$\ }}
\put(10,31){\makebox(0,0)[b]{$x_1$\ }}
\put(11,19){\makebox(0,0)[lt]{$y_1$\ }}
\put(11,09){\makebox(0,0)[lt]{$z_1$\ }}
\put(10,-1){\makebox(0,0)[t]{$x_1$\ }}
\put(20,31){\makebox(0,0)[b]{$x_2$\ }}
\put(21,19){\makebox(0,0)[lt]{$y_2$\ }}
\put(21,09){\makebox(0,0)[lt]{$z_2$\ }}
\put(20,-1){\makebox(0,0)[t]{$x_2$\ }}
\put(30,30){\makebox(0,0)[l]{\ $a$}}
\put(30,20){\makebox(0,0)[l]{\ $c$}}
\put(30,10){\makebox(0,0)[l]{\ $b$}}
\put(30,00){\makebox(0,0)[l]{\ $a$}}
\put(15,-10){\makebox(0,0)[b]{(a)}}
\end{picture}
}


\put(45,0){
\begin{picture}(40,45)(-5,-10)
\thicklines
\multiput(0,0)(10,0){4}{\circle*{1}}
\multiput(0,10)(10,0){4}{\circle*{1}}
\multiput(0,20)(10,0){4}{\circle*{1}}
\multiput(0,30)(10,0){4}{\circle*{1}}
\multiput(0,0)(0,10){4}{\line(1,0){30}}
\multiput(0,0)(10,0){4}{\line(0,1){30}}
\put(00,20){\line(1,1){10}}
\put(10,20){\line(1,1){10}}
\put(00,10){\line(1,1){10}}
\put(10,10){\line(1,1){10}}
\put(20,10){\line(1,1){10}}
\put(10,00){\line(1,1){10}}
\put(20,00){\line(1,1){10}}
\put(0,30){\makebox(0,0)[r]{$a$\ }}
\put(0,20){\makebox(0,0)[r]{$b$\ }}
\put(0,10){\makebox(0,0)[r]{$c$\ }}
\put(0,00){\makebox(0,0)[r]{$a$\ }}
\put(10,31){\makebox(0,0)[b]{$x_1$\ }}
\put(11,19){\makebox(0,0)[lt]{$y_1$\ }}
\put(11,09){\makebox(0,0)[lt]{$z_1$\ }}
\put(10,-1){\makebox(0,0)[t]{$x_1$\ }}
\put(20,31){\makebox(0,0)[b]{$x_2$\ }}
\put(21,19){\makebox(0,0)[lt]{$y_2$\ }}
\put(21,09){\makebox(0,0)[lt]{$z_2$\ }}
\put(20,-1){\makebox(0,0)[t]{$x_2$\ }}
\put(30,30){\makebox(0,0)[l]{\ $a$}}
\put(30,20){\makebox(0,0)[l]{\ $c$}}
\put(30,10){\makebox(0,0)[l]{\ $b$}}
\put(30,00){\makebox(0,0)[l]{\ $a$}}
\put(15,-10){\makebox(0,0)[b]{(b)}}
\end{picture}
}

\end{picture}

\caption{Partial structure PS1.1 and required edges.}
\label{PS1.1}
\end{center}
\end{figure}
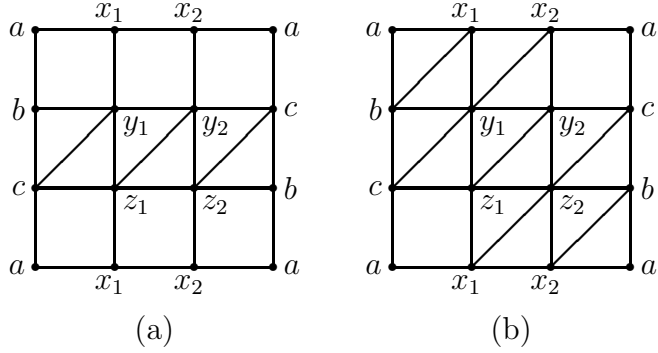

Consider the partial structure PS1.1 (Fig.~\ref{PS1.1}a).
The edge $by_1$ is a diagonal of the polygon $bax_1y_1z_1c$. 
By Lemma 3 of~\cite{LN} there must be a path of length 2 outside 
this polygon connecting $b$ and $y_1$.
Thus $y_1x_2$ and $x_2b$ must be edges of PS1.1. 
Likewise, $bx_1$ and $x_1z_2$ must be edges of PS1.1 because of
the diagonal $z_2b$ in the polygon $z_2y_2cbax_2$. 
We then have the partial structure in Fig.~\ref{PS1.1}b.  
If $y_2a$ is a diagonal of the polygon $y_2x_2ac$ in Fig.~\ref{PS1.1}b,
then $az_1$ must also be an edge of PS1.1.
Also, if $az_1$ is a diagonal of the polygon $acz_1x_1$ in 
Fig.~\ref{PS1.1}b,
then $y_2a$ must also be an edge of PS1.1.
So the two possible completed structures from PS1.1 are Kh14 and Kh15.

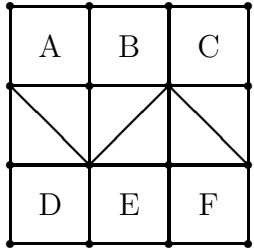
\begin{figure}
\begin{center}
\setlength{\unitlength}{3pt}

\begin{picture}(30,30)(0,0)
\thicklines
\multiput(0,0)(0,10){4}{\line(1,0){30}}
\multiput(0,0)(10,0){4}{\line(0,1){30}}
\multiput(0,0)(10,0){4}{\circle*{1}}
\multiput(0,10)(10,0){4}{\circle*{1}}
\multiput(0,20)(10,0){4}{\circle*{1}}
\multiput(0,30)(10,0){4}{\circle*{1}}
\put(0,20){\line(1,-1){10}}
\put(10,10){\line(1,1){10}}
\put(20,20){\line(1,-1){10}}
\put(5,25){\makebox(0,0){A}}
\put(15,25){\makebox(0,0){B}}
\put(25,25){\makebox(0,0){C}}
\put(5,5){\makebox(0,0){D}}
\put(15,5){\makebox(0,0){E}}
\put(25,5){\makebox(0,0){F}}

\end{picture}

\caption{Partial structure PS1.2.}
\label{PS1.2}
\end{center}
\end{figure}

We now determine the number of nonequivalent ways to complete the 
partial structure PS1.2, thus avoiding checking all $2^6$ configurations.
We examine two operations on PS1.2 which maintain its basic structure 
but which permute the quadrilateral regions labeled with the 
letters shown in Fig.~\ref{PS1.2}.
Reading from left to right in Fig.~\ref{PS1.2} and recalling that 
the left end of $R_{02}$ is identified with the right end of $R_{01}$, 
the order of the letter is ABCDEF. 
If PS1.2 is reflected both vertically and
horizontally, then it is still a type PS1.2 partial structure. 
The order of the letters is now FEDCBA and the permutation is 
$(AF)(BE)(CD)$.
If the left two columns of PS1.2 are removed,
reflected vertically, and pasted on the right, then the result
is still a type PS1.2 partial structure.  
The order of the letters is now CDEFAB and the permutation is 
$(ACE)(BDF)$.
These two operations define a permutation group 

\begin{eqnarray*}
\{(A)(B)(C)(D)(E)(F), (ACE)(BDF), (AEC)(BFD), \\
(AF)(BE)(CD), (AB)(CF)(DE), (AD)(BC)(EF)\}
\end{eqnarray*}

\noindent
acting on $\{A,B,C,D,E,F\}$.  
We color each element of $\{A,B,C,D,E,F\}$ with one of two colors 0 and 1.
Select a letter from $\{A,B,C,D,E,F\}$.
This letter is the label of one quadrilateral region in $R_{01}$
or $R_{02}$ which is adjacent to one quadrilateral region in $R_{12}$ 
which has a fixed diagonal in PS1.2.  
If the diagonals in these two quadrilateral regions
share a vertex (they are ``perpendicular''), 
then the letter is given the color 1.
If the diagonals in these two quadrilateral regions do not 
share a vertex (they are ``parallel''), 
then the letter is given the color 0.
The coloring has been defined so that the colors do not change 
under the group operations described above.
It can be shown that the number of equivalence classes of colorings 
is 16.


\begin{table}
\caption{Complete structures obtained from PS1.2}
\label{table1}
\vspace{.1in}
\begin{tabular}{r|l}
\hline
Coloring & Complete structure \\ \hline  
000000 & equivalent to Kh14  \\
000001 & has a contractible edge \\
000011 & equivalent to Kh15  \\
000101 & has a contractible edge \\
000110 & equivalent to Kh24 \\
001001 & has two contractible edges \\
111000 & Kh19 \\
001101 & has a contractible edge \\
101010 & Kh23 \\
100101 & Kh22 \\
101101 & Kh24 \\
111001 & Kh20 \\
111010 & Kh21 \\
110011 & Kh17 \\
110111 & Kh18 \\
111111 & Kh16 \\
\hline  
\end{tabular}

\end{table}

Table~\ref{table1} lists one coloring from each equivalence class 
along with
the resulting complete structure obtained from PS1.2.
From this table we see that 
the partial structure PS1.2 and thus PS1 produces the irreducible 
triangulations Kh14 through Kh24 and no others.

\section{Remarks}
\label{remarks}

In this section we reconsider three results appearing in \cite{BNN,LN,LN2}
whose proofs are based 
on the list of irreducible triangulations of the Klein bottle. Their proofs 
require that the irreducible triangulations of the Klein bottle 
have certain properties.  
Since, as we will observe, the additional triangulations, Kh22 
through Kh25, also have these properties, these results remain valid.

Firstly, Theorem 10 of~\cite{LN} states that
{\em if an irreducible triangulation of the Klein bottle can be 
embedded in the torus, then it is equivalent to \textnormal{Kh1}}.
In its proof the 
set of all the irreducible triangulations of the Klein bottle are 
partitioned into four subsets.  
The structure of the triangulations in each subset is 
considered in turn.  
Since each of the additional triangulations, Kh22 through Kh25,
can be placed in one of these four subsets,
the proof can be easily modified.  This result is then used by
Lawrencenko and Negami 
in~\cite{LN2} to construct all graphs which are triangulations of both
the torus and the Klein bottle.

Secondly, by checking Kh22 through Kh25 we see, as was observed in~\cite{LN},
that every irreducible 
triangulation of the Klein bottle still includes:

\begin{itemize}
\item a disjoint pair of longitudes and a meridian which crosses 
each of the longitudes only once,
\item a meridian and an equator which cross each other at precisely 
two vertices,
\item a Hamilton cycle which is trivial on the Klein bottle,
\item a Hamilton cycle which is a meridian,
\item a Hamilton cycle which is a longitude, and
\item a Hamilton cycle which is an equator.
\end{itemize}

Thirdly, Theorem 12 of~\cite{LN} states that
{\em a triangulation of the Klein bottle includes two disjoint meridians 
if and only if it does not include an equator of length~3}.
Its proof uses the
fact that every irreducible triangulation of handle type includes
two disjoint meridians.  
This fact is true for Kh22 through Kh25; hence the proof needs no change.
Theorem 12 of~\cite{LN} is used by
Brunet, Nakamoto, and Negami
in \cite{BNN} to prove that every 5-connected
triangulation of the Klein bottle has a Hamilton cycle which is 
contractible.

\begin{ack}
The author is grateful to the referees for their helpful comments.
\end{ack}

\bibliographystyle{amsplain}
\providecommand{\bysame}{\leavevmode\hbox to3em{\hrulefill}\thinspace}
\providecommand{\MR}{\relax\ifhmode\unskip\space\fi MR }
\providecommand{\MRhref}[2]{%
  \href{http://www.ams.org/mathscinet-getitem?mr=#1}{#2}
}
\providecommand{\href}[2]{#2}

\end{document}